\def\vac{|{\rm vac}\rangle}
\def\R{{\bf R}}
\def\sh{{\rm sh}\,}

\def\bea{\begin{eqnarray}}
\def\ena{\end{eqnarray}}

\newtheorem{prop}{Proposition}
\def\C{{\bf C}}
\def\Z{{\bf Z}}

\documentstyle{article}
\begin{document}
\begin{center}
\vspace*{10mm}
{\Large \bf The integral formula for the solutions of
the quantum Knizhnik-Zamolodchikov equation
associated with $U_q(\widehat{sl}_n)$ for $|q|=1$\\[8mm]}
{Tetsuji Miwa and Yoshihiro Takeyama\\[8mm]
Research Institute for Mathematical Sciences,\\[0mm]
Kyoto University, Kyoto 606, Japan\\[10mm]
{\it Dedicated to the memory of Moshe Flato.}
}
\end{center}
\begin{abstract}
\noindent
We write the integral formula of Tarasov-Varchenko type
for the solutions to the quantum Knizhnik-Zamolodchikov
equation associated with a tensor product of
the vector representations of $sl_n$.
We consider the case where the deformation parameter $q$
satisfies $|q|=1$. We use the bosonization of the type II vertex operators
in order to find the hypergeometric pairing in this setting.
\end{abstract}

\section{Introduction}
In this paper we construct a family of solutions to the quantum
Knizhnik-Zamolodchikov equation \cite{FR}
associated with the vector representation of
$sl_n$. To be more precise, we consider the case
\bea
\lambda&=&{4\pi\over n}
\label{LAM}
\ena
(called level $0$) of the following difference equation:
\bea
&&f(\beta_1,\ldots,\beta_r-\lambda i,\ldots,\beta_N)
=R_{r,r-1}(\beta_r-\lambda i,\beta_{r-1})
\cdots R_{r,1}(\beta_r-\lambda i,\beta_1)\nonumber\\
&&\phantom{-----}\times
D_rR_{r,N}(\beta_r,\beta_N)\cdots R_{r,r+1}(\beta_r,\beta_{r+1})
f(\beta_1,\ldots,\beta_N).\label{QKZ}
\ena
In this equation, the unknown function $f(\beta_1,\ldots,\beta_N)$
$(\beta_1,\ldots,\beta_N\in\C)$ takes its value
in $\underbrace{V\otimes\cdots\otimes V}_N$ where
\bea
V&=&\oplus_{j=0}^{n-1}\C v_j\label{V}
\ena
is the vector representation of $sl_n$.
The matrix $R_{rs}(\beta_r,\beta_s)$ acts on the $r$-th and $s$-th components
of the tensor product. It is the $R$-matrix for $U_q(\widehat{sl}_n)$.
The explicit formula of the $R$-matrix is given by (\ref{R}). 
In this paper, we restrict to the case where
\bea
q&=&e^{-{2\pi^2i\over\rho n}}\quad(\rho\in\R_{>0}).
\label{Q}
\ena
The matrix $D_r$ is a diagonal matrix acting on the $r$-th component.
In this paper, we restrict to the case
\bea
D_r&=&1.
\label{DIA}
\ena

For $n=2$, Smirnov \cite{S} constructed a family of solutions to (\ref{QKZ})
which he identified with the form factors of the quantum sine-Gordon model.
Later, Lukyanov \cite{L} gave a construction of the form factors by using
the bosonic vertex operators. In the $sl_n$ case,
we adapt Lukyanov's approach with the modification
given in \cite{JKM}. In this paper, however, we only construct the
so-called type II vertex operators. This construction gives us the
hypergeometric pairing for two functions $w$ and $W$
in the terminology of Tarasov and Varchenko \cite{TV1,TV2,MT}:
\bea
I(w,W)&=&\prod_{j=1}^{n-1}\prod_{1\leq r\leq \nu_j}
\int{d\gamma_{j,r}\over2\pi i}K(\{\gamma_{j,r}\})
w(\{\gamma_{j,r}\})W(\{\gamma_{j,r}\}).\label{HGP}
\ena

In this formula, the dependence on $\beta_1,\ldots,\beta_N$
in the right hand side is implicit in the kernel $K$ as well as $w$ and $W$.
We write the integral formula of Tarasov-Varchenko type for the solutions
to (\ref{QKZ}) with $\lambda={4\pi\over n}$ by using this pairing:
\bea
f(\beta_1,\ldots,\beta_N)
&=&\sum_{j_1,\ldots,j_r}I(w_{j_1,\ldots,j_r},W)v_{j_1}\otimes\cdots\otimes
v_{j_N}.
\label{FUN}
\ena
We prove the convergence of the integral (\ref{HGP})
in the positive Weyl chamber, and show that the equation (\ref{QKZ})
(with $\lambda={4\pi\over n}$) is indeed satisfied.

Finally we remark that the construction for $sl_n$ case
with generic highest weights is considered in \cite{TV}
using the Jackson integral, and in \cite{TV3}
in the setting of the hypergeometric pairing.
(see \cite{Matsuo, VS1, VS2, Kuroki} for the conformal case, i.e., $q=1$).

\bigskip
{\bf Acknowledgement}

The authors thank A. Nakayashiki, F. Smirnov and V Tarasov
for valuable discussions. One of the authors (T.M.) thanks the organizers
of the SIDE III for invitation and giving the opportunity for writing
this paper.  He also wishes to express his sorrow at the sudden death
of Moshe Flato. He shares the warmest memories of the days
with Moshe in Kyoto and Paris.

\def\shr{\sh{\pi\over\rho}}
\section{Operator construction}
In this section, we construct the vertex operators
\bea
\Psi^*_j(\beta)\quad(j=0,1,\ldots,n-1;\beta\in\R)
\ena
satisfying the commutation relations with the $R$-matrix for the
quantum affine algebra $U_q(\widehat{sl}_n)$, i.e.,
\bea
\Psi^*_j(\beta_1)\Psi^*_j(\beta_2)
&=&s(\beta_1-\beta_2)\Psi^*_j(\beta_2)\Psi^*_j(\beta_1),\label{CRA}
\ena
and for $j\not=k$
\bea
\Psi^*_j(\beta_1)\Psi^*_k(\beta_2)
&=&s(\beta_1-\beta_2)
\{\bar R(\beta_1,\beta_2)^{jk}_{jk}\Psi^*_k(\beta_2)\Psi^*_j(\beta_1)\nonumber\\
&+&\bar R(\beta_1,\beta_2)^{kj}_{jk}\Psi^*_j(\beta_2)\Psi^*_k(\beta_1)\}\label{CRB}
\ena
where
\bea
s(\beta)&=&
{S_2(-i\beta|\rho,2\pi)S_2(i\beta+{2(n-1)\pi\over n}|\rho,2\pi)
\over
S_2(i\beta|\rho,2\pi)S_2(-i\beta+{2(n-1)\pi\over n}|\rho,2\pi)},\\
\bar R(\beta_1,\beta_2)^{jk}_{jk}
&=&-{\sh{\pi\over\rho}(\beta_1-\beta_2)
\over\sh{\pi\over\rho}(\beta_1-\beta_2-{2\pi i\over n})}\quad(j\not=k),
\label{WB}\\
\bar R(\beta_1,\beta_2)^{kj}_{jk}
&=&\cases{
\displaystyle
-{e^{{\pi\over\rho}(\beta_1-\beta_2)}\sh{2\pi^2i\over\rho n}
\over\sh{\pi\over\rho}(\beta_1-\beta_2-{2\pi i\over n})}\quad(j>k);\cr
\displaystyle
-{e^{{\pi\over\rho}(\beta_2-\beta_1)}\sh{2\pi^2i\over\rho n}
\over\sh{\pi\over\rho}(\beta_1-\beta_2-{2\pi i\over n})}\quad(j<k).\cr}
\label{WC}
\ena
We define $\bar R(\beta_1,\beta_2)\in{\rm End}(V\otimes V)$ by

\bea
\bar R(\beta_1,\beta_2)v_j\otimes v_k&=&
\sum_{j'k'}\bar R(\beta_1,\beta_2)^{jk}_{j'k'}v_{j'}\otimes v_{k'},\label{R}
\ena
where $\bar R(\beta_1,\beta_2)^{jk}_{j'k'}=0$ except for
(\ref{WB}), (\ref{WC}) and $\bar R(\beta_1,\beta_2)^{jj}_{jj}=1$.

We refer the reader to \cite{JM} for the double sine function
$S_2(x|\omega_1,\omega_2)$. The deformation parameters $q$ and $\rho$ are
identified by the relation (\ref{Q}).
The construction of the vertex operators forces us to choose the normalization
factor $s(\beta_1-\beta_2)$.

The bosonization of the level $1$ highest weight representations for
$U_q(\widehat{sl}_n)$ and the vertex operators acting on them are given
in \cite{FJ} and \cite{K}, respectively. However, the present case is
not connected to those works in which $|q|<1$. The case $|q|=1$ is
obtained in the limit of the elliptic case considered in \cite{LP, AJMP}.
We adapt the construction in \cite{AJMP} by using the method
developed in \cite{L, JM, JKM}.

Let $a_j(t)$ $(1\leq j\leq n-1;t\in\R)$ be the free bose fields
satisfying the commutation relations
\bea
[a_j(t),a_k(t')]=A_{jk}(t)\delta(t+t'),\quad
A_{jk}(t)=-{1\over t}{\sh{a_{jk}\pi t\over n}\sh({\rho\over2}+{\pi\over n})t
\over\sh{\pi t\over n}\sh{\rho t\over2}}.
\ena
Here $(a_{jk})_{1\leq j,k\leq n-1}$ is the Cartan matrix of type $A_n$.

We consider the Fock space ${\cal F}$ generated by the vacuum vector $\vac$
satisfying
\bea
a_j(t)\vac=0\hbox{ if $t>0$}.
\ena

Set
\bea
a^*_1(t)&=&-\sum^{n-1}_{j=1}a_j(t){\sh{(n-j)\pi t\over n}\over\sh\pi t}.
\ena
We have
\bea
[a^*_1(t),a_j(t')]&=&\delta_{1j}{1\over t}{\sh({\rho\over2}+{\pi\over n})t
\over\sh{\rho t\over2}}\delta(t+t').
\ena

We introduce the currents
\bea
V_j(\beta)&=&:\exp\left(\int^\infty_{-\infty}a_j(t)e^{i\beta t}dt\right):
\quad(1\leq j\leq n-1).
\ena
They satisfy the following commutation relations.
\bea
V_j(\beta_1)V_k(\beta_2)&=&V_k(\beta_2)V_j(\beta_1)\quad(|j-k|\geq2),
\label{CR1}\\
V_j(\beta_1)V_{j-1}(\beta_2)
&=&{\sh{\pi\over\rho}(\beta_1-\beta_2+{\pi i\over n})
\over\sh{\pi\over\rho}(\beta_2-\beta_1+{\pi i\over n})}
V_{j-1}(\beta_2)V_j(\beta_1)\quad(2\leq j\leq n-1),\nonumber\\\label{CR2}\\
V_j(\beta_1)V_j(\beta_2)
&=&-{\sh{\pi\over\rho}(\beta_1-\beta_2-{2\pi i\over n})
\over\sh{\pi\over\rho}(\beta_2-\beta_1-{2\pi i\over n})}
V_j(\beta_2)V_j(\beta_1)\quad(1\leq j\leq n-1).\nonumber\\\label{CR3}
\ena

Now we define the vertex operators.
\bea
\Psi^*_j(\beta)&=&
\prod_{1\leq k\leq j}\int_{-\infty}^\infty{d\alpha_k\over2\pi i}
{e^{{\pi\over\rho}(\beta-\alpha_j)}
V_0(\beta)V_1(\alpha_1)\cdots V_j(\alpha_j)
\over\prod_{k=1}^j\sh{\pi\over\rho}(\alpha_{k-1}-\alpha_k+{\pi i\over n})}.
\label{VO}
\ena
Here we use $\alpha_0=\beta$ and
\bea
V_0(\beta)&=&
:\exp\left(\int^\infty_{-\infty}a^*_1(t)e^{i\beta t}dt\right):.
\ena
The relations (\ref{CR1}) and (\ref{CR2}) are also valid for $V_0(\beta)$.
However, the relation (\ref{CR3}) is modified to
\bea
V_0(\beta_1)V_0(\beta_2)
&=&s(\beta_1-\beta_2)V_0(\beta_2)V_0(\beta_1).
\ena
We give the proof of (\ref{CRA}) and (\ref{CRB}) later, which is
the same thing with the proof of Lemma \ref{lem1} in Section 3.

For a parameter $\lambda\in\R_{>0}$
we define the $\lambda$-expectation value by
\bea
\langle a_j(t)a_k(t')\rangle_\lambda=
{e^{\lambda t}\over e^{\lambda t}-1}A_{jk}(t)\delta(t+t').
\ena
If $\lambda=\infty$ this reduces to the vacuum expectation value.
The $\lambda$-expectation values of products of vertex operators
can be calculated by using Wick's theorem. We recall the following
formula (see \cite{JKM}) for the two point function:
\bea
\langle :e^{\int_{-\infty}^\infty a(t)e^{i\beta_1t}}:
:e^{\int_{-\infty}^\infty b(t)e^{i\beta_2t}}:\rangle_\lambda
&=&{\rm exp}\Bigl(\int_0^\infty A(t)
{{\rm ch}(i(\beta_1-\beta_2)+{\lambda\over2})t\over{\rm sh}
{\lambda t\over2}}dt\Bigr),\nonumber\\\label{WICK}
\ena
where $a(t)$ and $b(t)$ are bosons satisfying $[a(t),b(t')]=A(t)\delta(t+t')$
and $A(t)=-A(-t)$. We also list the two point
functions of $V_j(\beta)$. In the formulas below, const. means
a constant independent of the spectral parameters.
\bea\
\langle V_0(\beta_1)V_0(\beta_2)\rangle_\lambda&=&E_\lambda(\beta_1-\beta_2),\\
E_\lambda(\beta)&=&{\rm const.}
{S_3(-i\beta)S_3(i\beta+\lambda)\over
S_3({2\pi\over n}+\rho-i\beta)S_3({2\pi\over n}+\rho+i\beta+\lambda)},\\
\hbox{where}&&S_3(\beta)=S_3(\beta|\rho,\lambda,2\pi).
\ena
We have
\bea
{E_\lambda(\beta)\over E_\lambda(-\beta)}=s(\beta),\quad
E_\lambda(\lambda i-\beta)=E_\lambda(\beta).
\ena
For $1\leq j\leq n-1$ we have
\bea
\langle V_j(\beta_1)V_{j-1}(\beta_2)\rangle_\lambda
&=&\langle V_{j-1}(\beta_1)V_j(\beta_2)\rangle_\lambda\\
&=&{\rm const.}\,\sh{\pi\over\rho}(\beta_1-\beta_2+{\pi i\over n})
\varphi(\beta_1-\beta_2),\\
\varphi(\beta)&=&{1\over S_2(i\beta-{\pi\over n}|\rho,\lambda)
S_2(-i\beta-{\pi\over n}|\rho,\lambda)},\label{VARPHI}\\
\langle V_j(\beta_1)V_j(\beta_2)\rangle_\lambda
&=&{\rm const.}\,
\psi(\beta_1-\beta_2)\sh{\pi\over\rho}(\beta_1-\beta_2-{2\pi i\over n})
h(\beta_1-\beta_2,\lambda),\nonumber\\\\
h(\beta,\lambda)&=&\sh{\pi\over\lambda}\beta
\sh{\pi\over\lambda}(\beta-{2\pi i\over n})
\sh{\pi\over\lambda}(\beta+{2\pi i\over n}),\\
\psi(\beta)&=&{1\over S_2(i\beta+{2\pi\over n}|\rho,\lambda)
S_2(-i\beta+{2\pi\over n}|\rho,\lambda)}.
\ena

The rest of the two point functions of $V_j(\beta)$'s are $1$.
These formulas are to be understood as analytic continuations
from the regions where the integrals (\ref{WICK}) are convergent.
Because of the existence of poles, if they are used as integrands,
we must pay a special attention to the choice of integration contours.

To see this point closely, let us compute
\bea
\langle \Psi^*_{j_N}(\beta_N)\cdots\Psi^*_{j_1}(\beta_1)\rangle
&=&E(\beta_1,\ldots,\beta_N)\prod_{j,r}\int_C{d\alpha_{j,r}\over2\pi i}
I_{j_1,\ldots,j_N}(\{\alpha_{j,r}\}).
\ena
We associated the integration variables $\alpha_{j,r}$
$(1\leq r\leq N;1\leq j\leq j_r)$ to the current $V_j(\alpha_{j,r})$
contained in the vertex operator $\Psi^*_{j_r}(\beta_r)$ (see (\ref{VO})).
We also set
\bea
\alpha_{0,r}&=&\beta_r,\\
{\cal N}_j&=&\{r;j_r\geq j\}.
\ena

The factor $E(\beta_1,\ldots,\beta_N)$ is given by the pair product
\bea
E(\beta_1,\ldots,\beta_N)&=&{\rm const.}\,\prod_{1\leq r<s\leq N}
E_\lambda(\beta_s-\beta_r)
\ena
The integrand consists of three parts,
\bea
I_{j_1,\ldots,j_N}(\beta_1,\ldots,\beta_N)
&=&K(\{\alpha_{j,r}\})g(\{\alpha_{j,r}\})W(\{\alpha_{j,r}\}),
\ena
where
\bea
K(\{\alpha_{j,r}\})&=&\prod_{j=1}^{n-1}\Bigl\{
\prod_{r\in{\cal N}_j\atop s\in{\cal N}_{j-1}}
\varphi(\alpha_{j,r}-\alpha_{j-1,s})
\prod_{1\leq r<s\leq N\atop r,s\in{\cal N}_j}
\psi(\alpha_{j,s}-\alpha_{j,r})\Bigr\},\nonumber\\
\\
g(\{\alpha_{j,r}\})&=&
\prod_{j=1}^{n-1}\Bigl\{
\prod_{r\in{\cal N}_j}e^{{\pi\over\rho}(\alpha_{j-1,r}-\alpha_{j,r})}
\prod_{{r\in{\cal N}_j\atop s\in{\cal N}_{j-1}}\atop r>s}
\shr(\alpha_{j,r}-\alpha_{j-1,s}+{\pi i\over n})\nonumber
\ena
\bea
&&\times
\prod_{{r\in{\cal N}_j\atop s\in{\cal N}_{j-1}}\atop r<s}
\shr(\alpha_{j-1,s}-\alpha_{j,r}+{\pi i\over n})
\prod_{r,s\in{\cal N}_j\atop r<s}
\shr(\alpha_{j,s}-\alpha_{j,r}-{2\pi i\over n})\Bigr\},\nonumber\\
\label{G}
\ena
\bea
W(\{\alpha_{j,r}\})&=&
\prod_{j=1}^{n-1}
\prod_{r,s\in{\cal N}_j\atop r<s}
h(\alpha_{j,s}-\alpha_{j,r},\lambda).
\ena

The poles of the integrand come from those of $K(\{\alpha_{j,r}\})$.
They are located at
\bea
\alpha_{j,r}-\alpha_{j-1,s}
&=&\pm({\pi i\over n}-\rho i\Z_{\geq0}-\lambda i\Z_{\geq0}),\\
\alpha_{j,s}-\alpha_{j,r}
&=&\pm({2\pi i\over n}+\rho i\Z_{\geq0}+\lambda i\Z_{\geq0}).
\ena
The contour $C_{j,r}$ for $\alpha_{j,r}$ is chosen so that

\medskip
the poles at
$\alpha_{j-1,s}+{\pi i\over n}-\rho i\Z_{\geq0}-\lambda i\Z_{\geq0}$
are below $C_{j,r}$,

\medskip
the poles at
$\alpha_{j-1,s}-{\pi i\over n}+\rho i\Z_{\geq0}+\lambda i\Z_{\geq0}$
are above $C_{j,r}$,

\medskip
the poles at
$\alpha_{j,s}+{2\pi i\over n}+\rho i\Z_{\geq0}+\lambda i\Z_{\geq0}$
are above $C_{j,r}$,

\medskip
the poles at
$\alpha_{j,s}-{2\pi i\over n}-\rho i\Z_{\geq0}-\lambda i\Z_{\geq0}$
are below $C_{j,r}$.

\medskip
These conditions are not compatible if all the poles are
really existent. However, in the actual situation for our matrix element
$\langle \Psi^*_{j_N}(\beta_N)\cdots\Psi^*_{j_1}(\beta_1)\rangle_\lambda$, this is not
the case: the poles at
\bea
\alpha_{j,s}-\alpha_{j,r}&=&\pm{2\pi i\over n}
\ena
are canceled by the zeros of $W(\{\alpha_{j,r}\})$.
It is easy to see that under this cancellation, a consistent choice of the
contours is possible.

In this paper we finish the operator theory at this point. We will not
discuss the type I vertex operators and their form factors (\cite{L, JKM}).
In the next section, however, following the idea of Tarasov and Varchenko
(see \cite{TV1,TV2,MT}), we will modify $W(\{\alpha_{j,r}\})$, and thereby
construct a family of solutions to the quantum
Knizhnik-Zamolodchiov equation in the special case where
\bea
\lambda&=&{4\pi\over n}.
\ena
Note that, in this case, the function $\psi(\beta)$ simplifies to
\bea
\psi(\beta)&=&{1\over2i\sh{n\over4}(\beta-{2\pi i\over n})}.
\label{PSI}
\ena

\def\bea{\begin{eqnarray}}
\def\ena{\end{eqnarray}}
\def\no{\nonumber}
\def\sh{\mathop{\rm sh}\nolimits}
\def\ch{\mathop{\rm ch}\nolimits}
\def\Re{\mathop{\rm Re}\nolimits}
\def\Im{\mathop{\rm Im}\nolimits}
\def\qKZ{\hbox{$q${\rm KZ }}}
\newtheorem{lem}[prop]{Lemma}
\newtheorem{theo}[prop]{Theorem}

\newenvironment{pf}{{\it Proof.}}

\setcounter{section}{2}

\section{Integral formula}

For non-negative integers $\nu_{1}, \cdots , \nu_{n-1}$ satisfying
\bea
N=\nu_0\ge\nu_{1} \ge \nu_{2} \ge \cdots \ge \nu_{n-1} \ge\nu_n=0,
\label{NU}
\ena
we denote by ${\cal Z}_{\nu_{1}, \cdots , \nu_{n-1}}$ 
the set of all the $N$-tuples
$J=(j_{1}, \cdots , j_{N}) \in \left( {\bf Z}_{\ge 0} \right)^{N}$ 
such that
\bea
\nu_{j}&=&\#{\cal N}_j.
\ena

Define $r_{j,m}$ $(0\leq j\leq n-1;1\leq m\leq\nu_j)$ as follows:
\bea
{\cal N}_j=\{r_{j,1}, \cdots , r_{j, \nu_{j}} \},
 \quad r_{j,1}< \cdots < r_{j,\nu_{j}}.
\ena
We have, in particular, $r_{0, m}=m$.
We make a correspondence between two sets of the integration variables:
\bea
\alpha^J_{j,r_{j,m}}&=&\gamma_{j,m}.
\ena

We set
\bea
w_{J}(\{\gamma_{j,m}\})&=&{\rm Skew}_{n-1}\circ\cdots \circ{\rm Skew}_{1}
g_{J}(\{\gamma_{j, m} \}), \\
g_{J}(\{\gamma_{j,m}\})&=&\prod_{j=1}^{n-1}\Bigl\{
\prod_{m=1}^{\nu_{j}}\Bigl\{
\prod_{r\in{\cal N}_j}e^{\frac{\pi }{\rho}
(\alpha^J_{j-1,r}-\alpha^J_{j,r})} \!\!\!\!\!\!
\no
\ena
\bea
\times \!\!\!\!
\prod_{m' \atop r_{j-1, m'}<r_{j, m}}\!\!\!
\sh{\frac{\pi }{\rho}(\gamma_{j, m}-\gamma_{j-1, m'}+\frac{\pi i}{n})} 
 \!\!\!\!\!\!\!
\prod_{m' \atop r_{j, m}<r_{j-1, m'}}\!\!\!
\sh{\frac{\pi}{\rho}(\gamma_{j, m}-\gamma_{j-1, m'}-\frac{\pi i}{n})}
\Bigr\}
\no
\ena
\bea
\times\prod_{1 \le m<m' \le \nu_{j}}
\sh{\frac{\pi}{\rho}(\gamma_{j, m'}-\gamma_{j, m}-\frac{2\pi i}{n})}\Bigr\},
\ena
where $\gamma_{0, m}=\beta_{m}$ and 
${\rm Skew}_{j}$ is the skew-symmetrization 
with respect to the variables $\{ \gamma_{j, m} \}_{m=1, \cdots , \nu_{j}}$:
\bea
{\rm Skew}_jX(\gamma_{j,1},\ldots,\gamma_{j,\nu_j})
&=&\sum_{\sigma\in S_{\nu_j}}{\rm sgn}(\sigma)
X(\gamma_{j,\sigma(1)},\ldots,\gamma_{j,\sigma(\nu_j)}).
\ena

Note that the above definition differs from (\ref{G}) by sign.
This change is necessary because we will choose different $W$ in the below.
Accordingly, we change the $R$ matrix:
\bea
R(\beta_1,\beta_2)^{j'k'}_{jk}=\cases{-\bar R(\beta_1,\beta_2)^{jk}_{jk}
&if $(j',k')=(j,k)$ and $j\not=k$;\cr
\bar R(\beta_1,\beta_2)^{j'k'}_{jk}&otherwise.\cr}
\ena
In the following, we abbreviate
\bea
w_{J}(\{\gamma_{j, r}\})=w_{J}(\beta_{k_{1}}, \cdots , \beta_{k_{N}})
\ena
when the dependence on $\{\gamma_{j, m}\}$ $(j\not=0)$ is irrelevant.

For the function $w_{J}$, the following equality holds.
\begin{lem} \label{lem1}
\bea
& & w_{(j_{1}, \cdots , j_{k+1}, j_{k}, \cdots , j_{N})}
(\beta_{1}, \cdots , \beta_{k+1}, \beta_{k}, \cdots , \beta_{N})={} \no \\
& & \sum_{j_{k}', j_{k+1}'} \!\!
     R(\beta_{k}, \beta_{k+1})_{j_{k}j_{k+1}}^{j_{k}'j_{k+1}'} 
     w_{(j_{1}, \cdots , j_{k}', j_{k+1}', \cdots , j_{N})}
     (\beta_{1}, \cdots , \beta_{k}, \beta_{k+1}, \cdots , \beta_{N}).
\label{eq3}
\ena
\end{lem}
\begin{pf}
It is enough to prove (\ref{eq3}) for $N=2$.
First, we prove
\bea
w_{(l, l)}(\beta_{1}, \beta_{2})=w_{(l, l)}(\beta_{2}, \beta_{1}), \quad l=0, 1, \cdots , n-1.
\label{eq4}
\ena
Note that
\bea
&&g_{(l,l)}(\beta_1,\beta_2)=e^{{\pi\over\rho}
(\beta_1+\beta_2-\gamma_{l,1}-\gamma_{l,2})}\prod_{j=1}^l
\Bigl\{\shr(\gamma_{j,2}-\gamma_{j-1,1}+{\pi i\over n})\nonumber\\
&&\quad\times\shr(\gamma_{j,1}-\gamma_{j-1,2}-{\pi i\over n})
\shr(\gamma_{j,2}-\gamma_{j,1}-{2\pi i\over n})\Bigr\}.
\ena
By using
\bea
&&\shr(\gamma_2-\beta_1+{\pi i\over n})\shr(\gamma_1-\beta_2-{\pi i\over n})
-(\beta_1\leftrightarrow\beta_2)
\nonumber\\
&&\quad=\shr(\beta_1-\beta_2)\shr(\gamma_2-\gamma_1+{2\pi i\over n})
\ena
repeatedly, we obtain
\bea
&&{\rm Skew}_{l-1}\circ\cdots\circ{\rm Skew}_0g_{(l,l)}(\beta_1,\beta_2)
=e^{{\pi i\over\rho}(\beta_1+\beta_2-\gamma_{l,1}-\gamma_{l,2})}
\shr(\beta_1-\beta_2)\nonumber\\
&&\quad\times\shr(\gamma_{l,2}-\gamma_{l,1}+{2\pi i\over n})
\shr(\gamma_{l,2}-\gamma_{l,1}-{2\pi i\over n})
\prod_{j=1}^{l-1}h(\gamma_{j,1}-\gamma_{j,2})
.\label{SKEW1}
\ena
This is symmetric with respect to $\gamma_{l,1},\gamma_{l,2}$, and
therefore, we have (\ref{eq4}).

Next, we show that
\bea
{\rm Skew}_l\circ\cdots\circ{\rm Skew}_1X_{(l+1,l)}(\beta_2,\beta_1)=0
\label{eq5}
\ena
where
\bea
&&X_{(l+1,l)}(\beta_2,\beta_1)=
g_{(l+1,l)}(\beta_2,\beta_1)\nonumber\\
&&\quad-R(\beta_1,\beta_2)^{l+1,l}_{l,l+1}g_{(l+1,l)}(\beta_1,\beta_2)
-R(\beta_1,\beta_2)^{l,l+1}_{l,l+1}g_{(l,l+1)}(\beta_1,\beta_2).
\ena
The proof for $g_{(l,l+1)}(\beta_2,\beta_1)$ is similar.
The proof for $g_{(l,k)}$ for general $j,k$ reduces to the case $k=l\pm1$.

For $l=0$, one can check (\ref{eq5}) directly.
For $l\ge1$ we proceed as follows.
Set $\gamma_{l,1}=\gamma_1,\gamma_{l,2}=\gamma_2,\gamma_{l+1,1}=\gamma$.
Noting that
\bea
g_{(l+1,l)}(\beta_1,\beta_2)&=&g_{(l,l)}(\beta_1,\beta_2)
e^{{\pi i\over\rho}(\gamma_1-\gamma)}\shr(\gamma-\gamma_2-{\pi i\over n}),\\
g_{(l,l+1)}(\beta_1,\beta_2)&=&g_{(l,l)}(\beta_1,\beta_2)
e^{{\pi i\over\rho}(\gamma_2-\gamma)}\shr(\gamma-\gamma_1+{\pi i\over n}),
\ena
we have
\bea
&&X_{l+1,l}(\beta_2,\beta_1)
=\Bigl({1\over2}e^{-{2\pi\gamma\over\rho}}
-e^{{\pi\over\rho}(\gamma_1-\gamma_2-{\pi i\over n})}\Bigr)
{\rm Skew}_0g_{(l,l)}(\beta_1,\beta_2)\nonumber\\
&&\quad+{2\shr(\beta_1-\beta_2)\shr(\gamma_1-\gamma_2-{2\pi i\over n})
e^{-{\pi^2 i\over\rho n}}
\over\shr(\beta_1-\beta_2-{2\pi i\over n})}g_{(l,l)}(\beta_1,\beta_2).
\ena
Note that
\bea
&&{\rm Skew}_1\circ\cdots,\circ{\rm Skew}_l
\Bigl\{\shr(\gamma_1-\gamma_2-{2\pi i\over n})g_{(l,l)}(\beta_1,\beta_2)\Bigr\}
\nonumber\\&&\quad=\shr(\beta_1-\beta_2-{2\pi i\over n})e^{{\pi i\over\rho}
(\beta_1+\beta_2-\gamma_1-\gamma_2)}\prod_{j=1}^lh(\gamma_{j,1}-\gamma_{j,2}).
\label{SKEW2}
\ena
Using (\ref{SKEW1}) and (\ref{SKEW2}), we obtain (\ref{eq5}).
$\Box$ \newline 
\end{pf}

Following \cite{TV1, TV2, MT}, we define the hypergeometric
pairing.
Set
\bea
{\cal F}_{\nu_{1},\cdots,\nu_{n-1}}^{(\rho)}
&=&\sum_{J\in{\cal Z}_{\nu_{1},\cdots,\nu_{n-1}}}\C w_{J}(\{\gamma_{j,m}\}).
\ena
In our setting, the hypergeometric pairing gives rise to the pairing
given by (\ref{HGP}) between $w\in{\cal F}_{\nu_{1},\cdots,\nu_{n-1}}^{(\rho)}$
and $W\in{\cal F}_{\nu_{1},\cdots,\nu_{n-1}}^{(\lambda)}$, where
the kernel $K$ is given by
\bea
\prod_{j=1}^{n-1}\Bigl\{
\prod_{m=1}^{\nu_j}\prod_{m'=1}^{\nu_{j-1}}
\varphi(\gamma_{j,m}-\gamma_{j-1,m'})
\prod_{1\leq m<m'\leq \nu_j}\psi(\gamma_{j,m}-\gamma_{j,m'})\Bigr\}.
\ena
In this paper, we restrict to the level $0$ case, i.e.,
$\lambda={4\pi\over n}$, where we have (\ref{PSI}). In this case,
it is convenient to move the $\psi$-part of the kernel to the space of
$W$. Therefore, we set
\bea
{\cal W}_{\nu_{1},\cdots,\nu_{n-1}}
&=&\sum_{J\in{\cal Z}_{\nu_{1},\cdots,\nu_{n-1}}}\C W_{J}(\{\gamma_{j,m}\})
\ena
where
\bea
W_{J}(\{\gamma_{j,m}\})&=&{\rm Skew}_{n-1}\circ\cdots\circ
{\rm Skew}_{1}G_{J}(\{\gamma_{j, m} \}),\\
G_J(\{\gamma_{j, m}\})&=&\prod_{j=1}^{n-1}\Bigl\{
\prod_{r\in{\cal N}_j}
e^{\frac{n}{4}(\alpha^J_{j-1,r}-\alpha^J_{j,r})}\no\\
\times\prod_{m=1}^{\nu_{j}}\!\!\!\!
\prod_{m' \atop r_{j-1, m'}<r_{j, m}}\!\!\!\!
&&\!\!\!\!\!\!\!\!\!\!\!\!\!\!\!\!
\sh{\frac{n}{4}(\gamma_{j, m}-\gamma_{j-1, m'}+\frac{\pi i}{n})}
\!\!\!\!\!\!\!\!
\prod_{m' \atop r_{j, m}<r_{j-1, m'}}\!\!\!\!
\sh{\frac{n}{4}(\gamma_{j, m}-\gamma_{j-1, m'}-\frac{\pi i}{n})}
\Bigr\}.\no\\
\ena
We define the pairing as follows.
For $w \in {\cal F}_{\nu_{1}, \cdots , \nu_{n-1}}^{(\rho )}$
and $W \in {\cal W}_{\nu_{1}, \cdots , \nu_{n-1}}$, we set
\bea
I(w, W)(\beta_{1}, \cdots , \beta_{N}) \!\!\!
&=&\!\!\!\left(\prod_{j=1}^{n-1}\prod_{m=1}^{\nu_{j}}\int_{C_{j}}d\gamma_{j,m}\right)
\prod_{j=1}^{n-1}\prod_{m=1}^{\nu_{j}}\prod_{m'=1}^{\nu_{j-1}}
\varphi(\gamma_{j,m}-\gamma_{j-1,m'})\no\\
&&\quad\times w(\{\gamma_{j,m}\})W(\{\gamma_{j,m}\}).
\label{INTW}
\ena
In the above formula, $\varphi$ is
given by (\ref{VARPHI}) with $\lambda =\frac{4\pi}{n}$. 
The contour $C_{j}$ for $\gamma_{j, m}$ is $(-\infty, \infty)$ 
except that the poles at
\bea
\gamma_{j \pm 1, m'}+\frac{\pi i}{n}-\rho{\bf Z}_{\ge 0}-\frac{4\pi}{n}{\bf Z}_{\ge 0}
\ena
are below $C_{j}$ and the poles at
\bea
\gamma_{j \pm 1, m'}-\frac{\pi i}{n}+\rho{\bf Z}_{\ge 0}+\frac{4\pi}{n}{\bf Z}_{\ge 0}
\ena
are above $C_{j}$.

The weight of the function $f(\beta_1,\ldots,\beta_N)$ given by (\ref{FUN}) is
\bea
\sum_{j=0}^{n-1}(\nu_j-\nu_{j+1})\varepsilon_j.
\ena
The positive Weyl chamber is given by
\bea
\nu_{j-1}+\nu_{j+1} \ge 2\nu_{j}, \quad \hbox{ for all $j=1,\cdots,n-1$}.
\label{convcon}
\ena
\begin{prop}
Suppose that $\nu_{1}, \cdots , \nu_{n-1}$ satisfy
$(\ref{NU})$ and $(\ref{convcon})$.
Then the integral $(\ref{INTW})$ is absolutely convergent.
\end{prop}
\begin{pf}
It is easy to see that
\bea
& & \left| w(\{\gamma_{j, m}\}) \right| \le{} \no \\
& &{\rm const.}\exp\Bigl\{ \frac{\pi}{\rho}\sum_{j=1}^{n-1}\left( 
                \sum_{m=1}^{\nu_{j}}\sum_{m'=1}^{\nu_{j-1}}
                       |\gamma_{j, m}-\gamma_{j-1, m'}|+ \!\!\!
                \sum_{1 \le m<m' \le \nu_{j}} \!\!\!
                       |\gamma_{j, m}-\gamma_{j, m'}| \right)\Bigr\}, \no \\
& & \left| W(\{\gamma_{j, m}\}) \right| \le{} \no \\
&&{\rm const.}\exp\Bigl\{\frac{n}{4}\sum_{j=1}^{n-1}
\left(\sum_{m=1}^{\nu_{j}}\sum_{m'=1}^{\nu_{j-1}}
|\gamma_{j,m}-\gamma_{j-1,m'}|\right)\Bigr\},\no
\ena
where const. is a constant independent of $\{\gamma_{j, m}\}$.

The asymptotic behaviour of $\varphi $ is as follows (see \cite{JM}):
\bea
\varphi(\beta)\sim
\exp\Bigl\{-\left(\frac{n}{4}+\frac{3\pi}{2\rho}\right)|\beta|\Bigr\},
\quad \beta \to \pm\infty .
\ena
Therefore, we have
\bea
|{\rm the \, integrand \, of} \, I(w, W) | \le 
{\rm const.}\exp\Bigl(-\frac{\pi}{2\rho}H_{\nu_{1},\cdots,\nu_{n-1}}
(\{\gamma_{j,m}\})\Bigr),
\ena
where
\bea
& & H_{\nu_{1}, \cdots , \nu_{n-1}}(\{ \gamma_{j, m}\})={} \no \\
& & \sum_{j=1}^{n-1}\left(
     \sum_{m=1}^{\nu_{j}}\sum_{m'=1}^{\nu_{j-1}}|\gamma_{j, m}-\gamma_{j-1, m'}|-
     2\sum_{1 \le m<m' \le \nu_{j}}|\gamma_{j, m}-\gamma_{j, m'}| \right). \label{eq1}
\ena

Set
\bea
\nu=\sum_{j=0}^{n-1}\nu_j.
\ena 
Let $\gamma_1,\ldots,\gamma_\nu$ be a renumbering of the variables
$\{\gamma_{j,m}\}_{0\leq j\leq n-1\atop1\leq m\leq \nu_j}$.
We consider the
integral in the region
\bea
\gamma_1<\cdots<\gamma_\nu.
\ena
Set
\bea
\gamma_{k_1}&=&{\rm min}(\beta_1,\ldots,\beta_N),\\
\gamma_{k_2}&=&{\rm max}(\beta_1,\ldots,\beta_N).
\ena
One can rewrite (\ref{eq1}) as
\bea
H_{\nu_{1},\cdots,\nu_{n-1}}(\{\gamma_{j,m}\})
&=&\sum_{l=1}^{\nu-1}M_l(\gamma_{l+1}-\gamma_l).
\ena
Let
\bea
\tilde\gamma_1<\cdots<\tilde\gamma_{\nu-N}
\ena
be the renumbering of the variables
$\{\gamma_{j,m}\}_{1\leq j\leq n-1\atop1\leq m\leq \nu_j}$.
We change the integration variables from $\tilde\gamma_k$ to $t_k$
by setting
\bea
t_k=\cases{\gamma_{k+1}-\gamma_k&if $1\leq k<k_1$;\cr
\gamma_{k+N}-\gamma_{k+N-1}&if $k_2-N<k\leq\nu-N$;\cr
\tilde\gamma_k&otherwise.\cr}
\ena

For the convergence, it is enough to show that
\bea
M_k>0\quad(1\leq k<k_1),
\label{M1}
\ena
and
\bea
M_{k+N-1}>0\quad(k_2-N<k\leq\nu-N).
\label{M2}
\ena
We will show (\ref{M1}). The proof of (\ref{M2}) is similar.

Set
\bea
x_j=\sharp\{m;\gamma_{j,m}\leq\gamma_k\}
\quad(0\leq j\leq n-1). 
\ena
Note that $x_0=0$ because $k<k_1$.
Note also that not all $x_j$ $(1\leq j\leq n-1)$ are zero.
We have
\bea
M_{k}=\sum_{j=1}^{n-1}(\nu_{j-1}+\nu_{j+1}-2\nu_{j})x_{j}+
2\left( \sum_{j=1}^{n-1}x_{j}^{2}-\sum_{j=1}^{n-2}x_{j}x_{j+1} \right).
\ena
The second term of the right hand side of the above formula is positive.
Therefore, $M_{k}>0$ if the condition (\ref{convcon}) is satisfied.
$\Box$ \newline
\end{pf}

The hypergeometric pairing has the following property.
\begin{lem} \label{lem2}
\bea
&&I\bigl(g_{\overline J}(\beta_N,\beta_1,\cdots,\beta_{N-1}),
W(\beta_{1},\cdots,\beta_{N-1},\beta_{N})\bigr)
\Big|_{\beta_N\rightarrow\beta_N-{4\pi i\over n}}\no\\
&&\quad=I\bigl(g_{J}(\beta_{1},\cdots,\beta_{N}),
W(\beta_{1},\cdots,\beta_{N})\bigr),
\label{eq11}
\ena
where $\overline{J}=(j_{N},j_{1},\cdots,j_{N-1})$ for $J=(j_{1},\cdots,j_{N})$, 
and the left hand side is understood as the analytic continuation of the integral.
\end{lem}
\begin{pf}
We prove (\ref{eq11}) in the form
\bea
&&I\bigl(g_{\overline J}(\beta_N,\beta_1,\cdots,\beta_{N-1}),
W(\beta_{1},\cdots,\beta_{N-1},\beta_{N})\bigr)\no\\
&&\quad=
I\bigl(g_{J}(\beta_{1},\cdots,\beta_{N}),W(\beta_{1},\cdots,\beta_{N})\bigr)
\Big|_{\beta_N\rightarrow\beta_N+{4\pi i\over n}}.
\label{ANOTHER}
\ena

Note that the integration variables $\gamma_{j,\nu_{j}}\,(j=1,\cdots,j_{N})$
are associated with the $N$-th component in the sense of the operator construction in
Section 2. In other words, we have
\bea
r_{j,\nu_j}&=&N.
\ena
Therefore, it is natural to shift the contours for these variables
at the same time when we change $\beta_N$. In fact, by this simultaneous shift,
no crossing of contours by poles occurs. To see this, it is enough to observe that
the function $g_J$ has zeros at $\gamma_{j\pm1,m}=\gamma_{j,\nu_j}+{\pi i\over n}$
if $r_{j\pm1,m}<N$.

We have
\bea
W(\{\gamma_{j,m}\}_{m\not=N})
\Big|_{\gamma_{j,\nu_j}\rightarrow\gamma_{j,\nu_j}+{4\pi i\over n}
\atop(0\leq j\leq j_N)}
&=&(-1)^{N+\nu_{j_{N}}+\nu_{j_{N}+1}-1}W(\{\gamma_{j, m}\}). \label{eq12.5}
\ena
Note also that (see \cite{JM})
\bea
\frac{\varphi(\beta+\frac{4\pi i}{n})}{\varphi(\beta )}=
-\frac{\sh{\frac{\pi}{\rho}(\beta-\frac{\pi i}{n})}}
{\sh{\frac{\pi}{\rho}(\beta+\frac{5\pi i}{n})}}.
\label{eq13}  
\ena

We shift the contours in the integral (\ref{INTW})
(with $w=g_J$), and then change the variables $\gamma_{j,r}$
$(1\leq j\leq j_N)$: first by
\bea
\gamma_{j,\nu_{j}}\to\gamma_{j,\nu_{j}}+\frac{4\pi i}{n},
\ena
and second by
\bea
\gamma_{j,r}\rightarrow\gamma_{j,r+1}\quad(1\leq r\leq \nu_j-1),
\quad\gamma_{j,\nu_j}\rightarrow\gamma_{j,1}.
\ena
Rewriting the integrand by using (\ref{eq12.5}) and (\ref{eq13}),
we get (\ref{ANOTHER}).
$\Box$ \newline
\end{pf}

{}From Lemma \ref{lem1} and Lemma \ref{lem2}, it is easy to show the
following theorem.
\begin{theo}
We asuume the condition $(\ref{convcon})$. 
For $W \in {\cal W}_{\nu_{1}, \cdots , \nu_{n-1}}$,
we set
\bea
f_{W}(\beta_{1}, \cdots , \beta_{N})=
\sum_{J \in {\cal Z}_{\nu_{1}, \cdots , \nu_{n-1}}}
I(w_{J}, W)(\beta_{1}, \cdots , \beta_{N})v_{J},
\ena
where $v_{J}=v_{j_{1}}\otimes \cdots \otimes v_{j_{N}}$ for $J=(j_{1}, \cdots , j_{N})$.
Then $f_{W}$ is a solution to $(\ref{QKZ})$ with the restriction $(\ref{LAM})$
and $(\ref{DIA})$.
\end{theo}

\end{document}